\documentclass[twoside,11pt]{amsart} 
\usepackage{amsfonts}

 \newtheorem{thm}{Theorem} 
 
\newtheorem{lemma}[thm]{Lemma}

\usepackage{graphicx}
\def\qq{\hbox{\sl I\kern-.5em Q \kern-.3em}} 
\def\pf {{\bf Proof:} \ }

\def\ms{\medskip}

\begin{document}

\title{Erratum:
Studying links via closed braids IV: composite links and split links}

\author{Joan S. Birman \& William W. Menasco}

\maketitle
\centerline{July 22, 2004}

\flushleft

The purpose of this erratum is to fill a  gap  in the proof of the `Composite Braid Theorem' in the
manuscript {\em Studying Links Via Closed Braids IV: Composite Links and Split Links},
Inventiones Math, {\bf 102} Fasc. 1 (1990), 115-139.  The statement of the theorem is unchanged. 
The gap occurs on page 135, lines $13^-$ to $11^-$, where we fail to consider the
case:
\begin{itemize}
\item [] $V_2 = 4, V_4 > 0, V_j=0$ if $j\not= 2,4,$ and all 4 vertices of valence 2 are bad. 
\end{itemize}

At the end of this Erratum we make some brief remarks on the literature, as it evolved during the 14 years
between the publication of \cite{SLVCB-IV} and the submission of this Erratum.

\begin{figure} [htpb]
\centerline{\includegraphics[scale=.48]{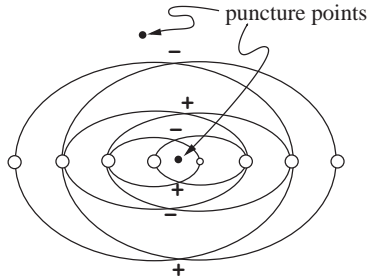}}
\caption{A tiling of the 2-sphere {\bf Y} with $V_2=V_4=4$, $V_j=0$ if $j\not= 2,4$, and all vertices of valence 2 are bad.}
\label{figure:example}
\end{figure}

See Figure \ref{figure:example} below for an example which illustrates that the missing
case can occur.  Ivan Dynnikov discovered it when he was working on his manuscript \cite{Dynnikov}, where he
established two theorems about arc presentations of links which are similar to the two theorems that we had
proved for closed braid presentations in \cite{SLVCB-IV}. His proof was a modification of our proof to new
geometry, and in the course of his work he realized there was a gap.  We thank him for pointing it out to us.

To handle the missing case, we need to change the definition of a {\em good vertex}, on page 135,
lines 23-24.  A vertex $v$ in the foliation of our $2$-sphere
is {\em good} if every non-singular leaf adjacent to $v$ has empty intersection with {\bf K}.  (As stated, we are
always assuming that {\bf K} intersects the foliation of  {\bf Y} away from singular leaves.)  With this change,
the definition of a {\em good tile} is no longer needed or used.  \ms

We now discuss how to handle the missing case, by showing that there is a new choice of disc
fibers for the fibration which, in effect, `moves one of the punctures out of a bad region', thereby
making it good. After that we will be able to reduce the complexity. The proof adapts  an
argument used at a later point in \cite{SLVCB-IV} to the situation that we overlooked.

\begin{figure} [htpb]
\centerline{\includegraphics[scale=.8]{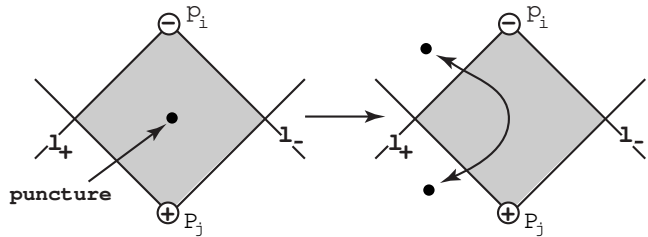}}
\caption{We illustrate the two ways that one of the puncture points
can be moved past the singular leaf. Both ways are realizable.}
\label{figure:disc region}
\end{figure}

A {\em disc region} $ R \subset {\bf Y}$
is a sub-disc of ${\bf Y}$ whose boundary is contained
in the union of two singular leaves and whose interior is foliated by a family of
$b$-arcs.  See Figure \ref{figure:disc region} below and Figure 6 of \cite{SLVCB-IV}.

\begin{lemma}
\label{lemma: change of fibration}
Let $R$ be a disc region in the tiling of {\bf
Y}, with $R \cap {\bf K} $  a single point, also $\partial R \subset l_+
\cup l_- $, where $l_\pm$ is a singular leaf whose associated singularity has
sign $\pm$.
Then there exists a new choice of disc fibers 
which leaves the the tiling of {\bf Y}  unchanged, but after the
change $R \cap {\bf K} =\emptyset$.  Moreover, if {\bf Y} is oriented, and if $R \cap
{\bf K} $ is a positive (resp. negative) puncture then the change of fibration corresponds to an
isotopy of ${\bf K}$ across $l_+$ (resp. $l_-$). As illustrated in Figure \ref{figure:disc region}, there
are two possible isotopies.
\end{lemma}

\pf  To simplify notation we consider a special case, assuming that:
\begin{itemize}
\item The negative singularity of $R$, the  $b$-arcs in $R$, and the
positive singularity of $R$ occur in that order in the fibration.
\item The puncture point  $R \cup {\bf K}$ is
positive.  
\end{itemize}
We will  establish that the order of
the positive  puncture and the positive singularity in the fibration can be
interchanged by re-choosing the disc fibers. 
The special case is in fact the general case because the other order and parity
possibilities can be realized by reversing the orientation of {\bf Y} and/or reversing
the orientation of the fibration. 

\begin{figure}[htpb]
\centerline{\includegraphics[scale=.85]{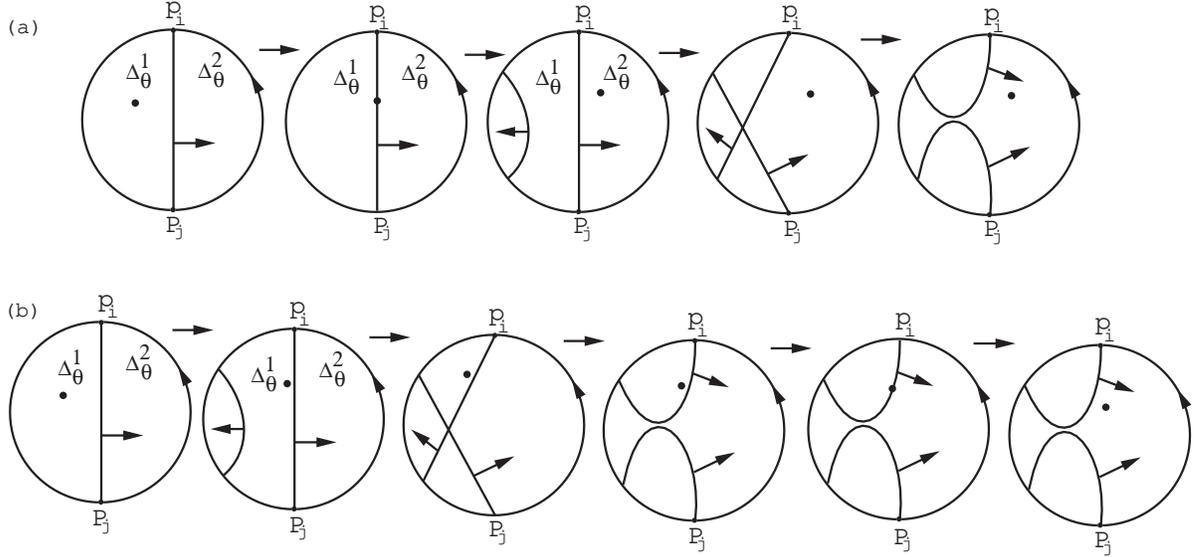}}
\caption{The two  $H_{\theta}$ sequences in the proof of Lemma 1}
\label{figure:sequences}
\end{figure}

\begin{figure}[htpb]
\centerline{\includegraphics[scale=.80]{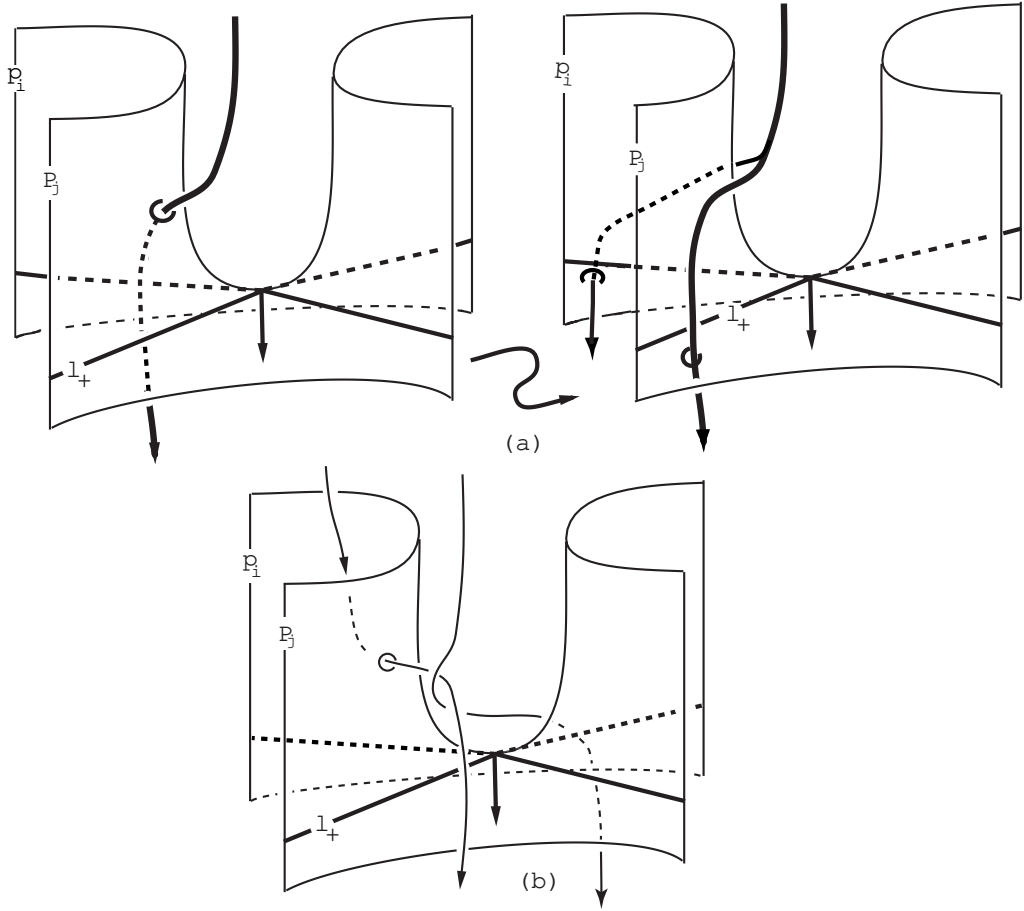}}
\caption{(a) The geometric realization of the change in fibration in Figure \ref{figure:disc region}. 
(b) We illustrate a possible obstruction to a change of fibration
when the parity of the puncture and singularity disagree.}
\label{figure:3-D}
\end{figure}

Referring to sketch (a) in Figure \ref{figure:sequences},  we depict the
$H_\theta$-sequence  from the instant just before the  positive puncture point to
the instant right after the positive singularity.  
Observe that the $b$-arcs split $H_\theta$ into two subdisc
$\Delta^1_\theta \cup \Delta^2_\theta = H_\theta$ where  $\Delta^1_\theta$  is on the
negative side of the $b$-arcs and $\Delta^2_\theta$ is on the positive of the
$b$-arcs. Thus the events---singularities and braiding---in the fibration that occur
in $\Delta^1_\theta$ are independent of the events that occur in $\Delta^2_\theta$.
This allows us to alter the $H_\theta$-sequence.   In
particular, prior to the instant when {\bf K} punctures {\bf Y} we can ``push''
forward the fibration on the $\Delta^1_\theta$ side until precisely after the
occurrence of the positive singularity. We then allow {\bf K} to puncture {\bf Y}.  
The ``dot'' in
$\Delta^1_\theta \cap {\bf K} \subset H_\theta$ that is used in (a) of Figure
\ref{figure:sequences} to create the puncture can be placed close to the
$b$-arcs of $R$ away from the other events as we do this ``push''.  We have a choice
of which side of the positive singularity we want the ``dot''; thus, when the
dot finally passes through {\bf Y}, it can be on either side of the
singularity. The sketches in (b) of Figure
\ref{figure:sequences} illustrates the resulting change of fibration.
Figure \ref{figure:3-D}(a) illustrates the (two possible) $3$-dimensional geometric realization of the
$H_\theta$-sequence of Figure \ref{figure:sequences}.  (For completeness we also
have Figure \ref{figure:3-D}(b) which illustrates a possible obstruction to a change of fibration
when the parity of the puncture and singularity disagree.)
All of the
other events are combinatorially unaltered by this new choice of disc
fibers, therefore the foliation of the {\bf Y} is unchanged. $||$

\begin{lemma}
\label{lemma: valence two vertices with punctures}  See Figure 13 of [1].
Let $link(p) = R_1 \cup R_2$ be the link neighborhood of
a valence two vertex $p$ in the tiling of the $2$-sphere ${\bf Y}$ which realizes the
connected sum where $R_1$ \& $R_2$ are disc regions.  Assume that
$|R_i
\cap {\bf K}| < 2, i=1,2$.  Then after the change of
fibration of Lemma 1 we can perform an exchange move on ${\bf K}$ and an isotopy of
{\bf Y} which eliminates $p$. 
\end{lemma}

\pf Suppose that
$R_1$ is a disc region in the link of vertex $p$ of valence two.  Referring to Figure 13 of [1],
we have two disc regions containing $b$-arcs, with $\{t_2, t_3 , t_4\}$ corresponding to
$R_1$ and $\{t_6, t_7 , t_8\}$ respectively corresponding  to
$R_2$.  Since $link(p)$ has both a positive and a negative singularity, the disc region
$R_1$ also has both positive and negative singular points in its boundary. By assumption
$R_1$ has at most one puncture point.  We apply Lemma \ref{lemma: change of fibration}  to
move the puncture point out of
$R_1$ (but not into $R_2$).  This is possible since as illustrated in Figure \ref{figure:disc region}
there are two possible change of fibration isotopies and we can choose which side of a
singular point the puncture passes by.  Similarly, for $R_2$.  After our application
of the lemma making $link(p)$ a good disc
we know we can apply an exchange move to eliminate the valence two vertex.  $||$

Lemmas 1 and 2 handle the missing case in the proof of the Composite Braid Theorem of [1]. $||$

In closing, we use this opportunity to correct an error in the signs of certain crossings  in Figure 22 of
\cite{SLVCB-IV}.   Each of the two crossings of the strand
that is labeled ``k'' (resp. the strand that is labeled ``m'') with the axis {\bf A} should be reversed. 

\

{\bf Remarks on the literature:} For the benefit of the reader, we note that during the 14 years between the
publication of \cite{SLVCB-IV} and the submission of this Erratum, there have been many papers published
which use the braid foliation techniques introduced in \cite{SLVCB-IV}, for example \cite{SLVCB-V},
\cite{SLVCB-VI}, \cite{Cromwell}, \cite{Dynnikov}, \cite{F-S} and \cite{Matsuda}.  For that reason, 8 years after
the manuscript
\cite{SLVCB-IV} appeared, the first author and E. Finkelstein wrote a review article on braid foliations
\cite{B-F}. It contains, in particular, careful proofs of the validity of the change in
fibration  and the exchange move, both of which first appeared in \cite{SLVCB-IV}. The interested reader may
find it helpful to consult \cite{B-F} if he/she is interested in learning the basic techniques in braid foliation
theory. 

In a different direction, readers have pointed out to us the
similarity between many of our techniques and those of D. Bennequin in \cite{Bennequin}. As one example, the
technique that we use in the proof of Lemma 3 of
\cite{SLVCB-IV} to reduce the valence of the two middle vertices in the sketches labled (c) in Figure 18
of
\cite{SLVCB-IV} is similar to the argument that is given in the proof of Theorem 6 of \cite{Bennequin}, in the
(unnumbered) Lemma that appears in the course of the proof of Theorem 6 of \cite{Bennequin}.  Those similarities
puzzled us greatly for many years, because Bennequin was studying the so-called characteristic foliation of a
surface in 3-space supplied with a contact struture, and not the braid foliation of the same surface in 3-space
supplied with a braid structure. When we first learned of such similarities we convinced ourselves that the
characteristic foliation and the braid foliation were genuinely different, and that most of the extensive
literature resulting from Bennequin's foundational work \cite{Bennequin} simply did not apply to braid
foliations.  After recent conversations with Ivan Dynnikov, we finally understood the explanation: the essential
hypotheses needed to prove particular overlapping technical lemmas are equally valid in the setting of braid
foliations and characteristic foliations.   Bennequin did not prove any of the theorems in
\cite{SLVCB-IV} or our subsequent papers based upon braid foliation techniques.  Moreover, it seems to us that
any applicaion of results worked out in either setting to the other setting would require a great deal of
careful checking, because the underlying setting, terminology, basic definitions and
references are all distinct.

\  

Joan S. Birman,  Dept of Math.,
Barnard College \& Columbia Univ., \\
2990 Broadway, NY, NY 10027\\  
e-mail jb@math.columbia.edu\\ 

\

 William W. Menasco, Dept of Math., SUNY at Buffalo, \\
Buffalo, NY 14222 \\
e-mail menasco@tait.math.buffalo.edu

\end{document}